\documentclass{article} 

\title{A partition theorem  for scattered  order types}

\author{P\'eter Komj\'{a}th\thanks{Research  partially 
                                   supported by Hungarian 
                                   National Research Grant T 032455.}\\
         Saharon Shelah\thanks{This research was supported by the 
         Israel  Science Foundation. Publication 796.}} 

\renewcommand\a{\alpha}
\renewcommand\b{\beta}
\newcommand\bb{\beta}
\newcommand\g{\gamma}
\newcommand\ga{\gamma}
\newcommand\dd{\delta}

\newcommand\la{\lambda}
\renewcommand\o{\omega}
\newcommand\oo{\omega}
\renewcommand\phi{\varphi}

\newcommand\fs{{\rm FS}}
\newcommand\tp{{\rm tp}}
\newcommand\qed{\hfill{\vbox{\hrule\hbox{\vrule\kern3pt
                \vbox{\kern6pt}\kern3pt\vrule}\hrule}}}

\begin{document}
\maketitle

\begin{abstract}
If $\phi$ is a scattered order type, $\mu$ a cardinal, then there exists a 
scattered order type $\psi$ such that $\psi \to [\phi]^{1}_{\mu,\aleph_0}$ 
holds. 
\end{abstract}

In this note we prove a Ramsey type statement on scattered order types. 
A trivial fact on ordinals implies the following statement. 
If $\mu$ is an infinite cardinal, then $\mu^+\to(\mu^+)^1_\mu$. 
It is less trivial but still easy to show that if $\phi$ is an 
order type, $\mu$ a cardinal then there is some order type $\psi$ 
that $\psi\to (\phi)^1_\mu$ holds. 
One can say that these results show that the classes of ordinals and 
order types are both Ramsey classes in the natural sense; given a 
target element and a cardinal for the number of colors, there is 
another element of the class, which, when colored with the required 
number of colors, always has a monocolored copy of the target. 
One can wonder which other classes have similar Ramsey properties. 
A natural, and well investigated, class in between is the class of 
{\em scattered order types}. 
For this class, the Ramsey property fails for the following well known and 
simple reason. 
There is some scattered order type $\psi$ that for every 
scattered $\phi$ one has 
$\phi \not\to [\psi]^{1}_{\o}$. 
See Lemma 1. 

In this paper we show that this is the most in the negative direction, 
that is, for every   scattered order type $\phi$ and cardinal $\mu$
there exists a 
scattered order type $\psi$ such that $\psi \to [\phi]^{1}_{\mu,\o}$ 
holds.

\medskip
\noindent{\bf Notation.} We use the standard axiomatic set theory notation. 
If $\phi$, $\psi$ are order types, then $\phi\leq \psi$ denotes 
that there is an order preserving embedding of $\phi$ into $\psi$, 
that is, every ordered set of order type $\psi$ has a subset of 
order type $\phi$. 
If $\phi$ is an order type, then $\phi^*$ denotes the reverse order type, 
that is, if $\phi$ is the order type of $(S,<)$, then $\phi^*$ 
is the order type of $(S,>)$. 
$\o$ is the ordinal of the set of natural numbers, $({\bf N},<)$. 
$\eta$ is the order type of the set of rational numbers, $({\bf Q},<)$.

If $\phi$, $\psi$ are order types, $\mu$ is a cardinal, 
$\phi \to (\psi)^{1}_{\mu}$ denotes the following statement. 
If $(S,<)$ is an ordered set of order type $\phi$ and 
$f:S\to\mu$ then for some $i<\mu$ the subset $f^{-1}(i)$ contains a subset 
of order type $\psi$. 
That is, if a set of order type $\phi$ is colored with $\mu$ colors, 
then there is a monochromatic $\psi$. 
If the statement does not hold, we cross the arrow, 
$\phi \not\to (\psi)^{1}_{\mu}$

If $\phi$, $\psi$ are order types, $\la$, $\mu$ cardinals, 
then 
$\phi \to [\psi]^{1}_{\la,\mu}$ denotes the following statement. 
If $(S,<)$ is an ordered set of order type $\phi$ and 
$f:S\to\la$ then there is a subset $X\subseteq \la$ of cardinality 
$\mu$ such that 
the set $\{x\in S:f(x)\in X\}$ contains a subset of order type $\psi$. 
Again, crossing the arrow denotes the negation of the statement; 
$\phi \not\to [\psi]^{1}_{\la,\mu}$. 
Notice that $\phi \to (\psi)^{1}_{\mu}$ is equivalent to 
$\phi \to [\psi]^{1}_{\mu,1}$.

If $\phi$, $\psi$ are order types, $\mu$ is a cardinal, 
$\phi \not\to [\psi]^{1}_{\mu}$ denotes the following statement. 
If $(S,<)$ is an ordered set of order type $\phi$ then there 
is a function $f:S\to\mu$ such that on every subset of $S$ 
of order type $\psi$, $f$ assumes every value. 
If the statement fails that is, we have a positive statement on 
all $f:S\to\mu$ function, then we do not cross the arrow; 
$\phi \to [\psi]^{1}_{\mu}$

\medskip
The order type $\phi$ is {\em scattered} iff $\eta\not\leq \phi$. 
Hausdorff proved that the class of scattered order types is exactly 
the smallest class containing $0$, $1$, and closed under 
well ordered and reversely well ordered sums 
(see \cite{erdhaj}, \cite{haus}, \cite{rosen}).

\medskip
\noindent{\bf Lemma 1.} 
         {\em If $S$ is an ordered set with the scattered order type 
         $\phi$ then there is some $f:S\to\o$ such that $f^{-1}(n)$ 
         has no subset of order type $\left(\o^*+\o\right)^n$. 
         Therefore, $\phi\not\to(\psi)^1_\o$ where 
         $\psi=1+\left(\o^*+\o\right)+\left(\o^*+\o\right)^2\cdots$.}

\medskip
\noindent{\bf Proof.}          
The second statement obviously follows from the first one. 
In order to prove the first statement, with Hausdorff 
characterization of scattered order types it suffices to show it 
for $(S,<)$ which is the well ordered sum of the 
ordered sets $\left\{(S_i,<):i<\a\right\}$ and we have the required function 
$f_i:S_i\to\o$ for every $i<\a$.

Define $f:S\to\o$ as follows. 
$f(x)=f_i(x)+1$ when $i<\a$ is the unique ordinal that 
$x\in S_i$. 
If we now have a set of order type $\left(\o^*+\o\right)^{n+1}$ in color 
$n+1$ then the $\o^*$ copies of $\left(\o^*+\o\right)^n$ in it  
left side must all but finitely many be in the same $S_i$, 
of color $n$, which contradicts the assumption on $f_i$. 
\qed

\medskip
Before proceeding to  our main theorem we need to show a technical result. 

\medskip
In what follows for an ordinal $\la$ we denote by 
$\fs(\la)$ the set of all finite decreasing sequences from $\la$, 
that is, an element ${\bf s}$ is of the form 
${\bf s}=s(0)s(1)\cdots s(n-1)$ with $\la>s(0)>s(1)>\cdots> s(n-1)$. 
Here $n=|{\bf s}|$ is the {\em length} of ${\bf s}$. 
The extension of the string ${\bf s}$ with one ordinal $\ga$ is denoted by 
${\bf s}\ga$. 
We therefore, identify finite subsets of $\la^+$ with 
decreasingly ordered strings. 

If $\a$ is an ordinal, then an {\em $\a$-tree} is a system of ordinals 
$\{x({\bf s}):{\bf s}\in \fs(\a)\}$ with the following properties: 

\[ 
      x({\bf s}\ga)< x({\bf s}\ga')<  x({\bf s})  \qquad 
      \mbox{for\ }\ga<\ga'<\min({\bf s}).
\] 

\medskip
\noindent{\bf Theorem 1.} 
         {\em Assume that $\a$ is an ordinal and $\mu$ is a cardinal. 
         Set $\la= \bigl(|\a|^{\mu^{\aleph_0}}\bigr)^+$. 
         Assume that $F:\fs(\la^+)\to\mu$ then there exist an $\a$-tree 
         $\{x({\bf s}):{\bf s}\in \fs(\a)\}$  and a function 
         $c:\oo\to\mu$ such that 
         \[ 
         F\biggl( x\bigl(s(0)\bigr), x\bigl(s(0)s(1)\bigr), \cdots,
         x\bigl(s(0)s(1)\cdots s(n)\bigr) \biggr)
%         x\bigl(s(0)s(1)\cdots s(n-2)\bigr), 
%         \dots, x\bigl(s(0)\bigr)
         =c(n)
         \] 
         holds for every element ${\bf s}=s(0)s(1)\cdots s(n)$ 
         of length $n+1$ of the 
         tree.}

\medskip
\noindent{\bf Proof.}          
We define, for every ${\bf s}\in \fs(\a)$ and for every function 
$c:\oo\to\mu$ a rank $r_c({\bf s})$ as follows. 
Assume that ${\bf s}=s(0)s(1)\cdots s(n-1)$. 
 $r_c({\bf s})=-1$ if for some $0\leq i < n$ we have 
$F(s(0)s(1)\cdots s(i))\neq c(i)$. 
Otherwise, we declare that $r_c({\bf s})\geq 0$. 
Then we define by induction on $\xi$ when $r_c({\bf s})\geq \xi$ 
holds; 
we set $r_c({\bf s})\geq \xi$ if and only if for every 
$\nu<\xi$ we have 
\[ 
      \la \leq \tp \bigl( \{ \ga < \min ({\bf s}): r_c({\bf s}\ga)\geq \nu
                          \}
                   \bigr). 
\] 
Naturally, $r_c({\bf s}) =  \xi$ holds if $r_c({\bf s})\geq \xi$  but 
$r_c({\bf s})\geq \xi+1$ is not true. 

Assume first that for some function $c:\oo\to\mu$ we have 
$r_c(\emptyset)\geq \a$. 
In this case we can select the $\a$-tree as required in the Theorem with 
the additional property that 
\[ 
r_c\biggl( x\bigl(s(0)\bigr), x\bigl(s(0)s(1)\bigr), \cdots,
         x\bigl(s(0)s(1)\cdots s(n)\bigr) \biggr)\geq s(n).       
\]

To show this we have to show that if we are given an ${\bf s}$ 
with $r_c({\bf s})\geq \bb$, then we can select the ordinals 
$\{x_\ga:\ga<\bb\}$ 
with $x_\ga<x_{\ga'}<\min({\bf s})$ for $\ga<\ga'<\bb$ and with 
$r_c({\bf s}x_\ga)\geq \ga$ for $\ga<\bb$. 
To this end, we let $\dd_\ga$ be the supremum of the first 
$\la$ ordinals $x$ with the property that $r_c({\bf s}x)\geq \ga$. 
Notice that $\dd_{\ga'}\leq\dd_\ga$ for $\ga'<\ga$ and the cofinality 
of is $\dd_\ga$ is $\la$. 
We are going to select by transfinite recursion the elements 
$x_\ga<\dd_\ga$ as required. 
At step $\ga$ we have the elements 
$\{x_{\ga'}:\ga'<\ga\}$ selected and as 
$\sup(\{x_{\ga'}:\ga'<\ga\})\leq \sup (\{\dd_{\ga'}:\ga'<\ga\})\leq\dd_\ga$ 
we have $\sup(\{x_{\ga'}:\ga'<\ga\})<\dd_\ga$ and so we can choose 
$x_\ga$.

Assume now that for every function $c:\oo\to\mu$ there holds 
$r_c(\emptyset)< \a$. 

In this case we construct by induction on $0\leq n<\oo$ the ordinals 
\[
      \{x(n,\ga,s):\ga<\la^+,s: k\to\la, k\leq n\},
\] 
the ordinals $d(n)<\mu$, and for every $c:\oo\to\mu$, 
the values $-1\leq\xi(n,c)<\a$ 
with the following properties 
\begin{equation} 
      x(n,\ga,s\tau)<x(n,\ga,s\tau')<x(n,\ga,s) 
      (1\leq |s|<n, \tau<\tau'<\min(s))
\end{equation}
\begin{equation} 
      \ga<x(n,\ga,s) 
\end{equation}
and finally, if $\ga<\la^+$, $s:n\to\la$, $1\leq k \leq n$, 
and we set  
      $y_i=x(n,\ga,s|i)$, then
\begin{equation} 
      F(y_0,\dots,y_k)=d(k)
\end{equation}    
and 
\begin{equation} 
      r_c(y_0,\dots,y_k)=\xi(k,c)
\end{equation}      
hold for every $c:\oo\to\mu$.

To start, we select $\la^+$ ordinals $x(0,\ga,\emptyset)$ ($\ga<\la^+$) such 
that  the value $F(x(0,\ga,\emptyset))$ is the same, let this be $d(0)$, and 
for every $c:\oo\to\mu$the value $r_c\bigl(x(0,\ga,\emptyset)\bigr)$ is the 
same, this will be $\xi(0,c)$. 
This is possible, by the pigeon hole principle, counting possibilities. 
%as we have less than $\la$ many possibilities for 
%the $c\mapsto r_c(x)$ functions, for $\la^+$ many of the ordinals $x$ 
%it is the same. 

Assume that we have the result for some value $n$ and we have the 
corresponding system 
$ \{x(n,\ga,s):\ga<\la^+,s:k\to\la, k\leq n\}$ with 
$\ga<x(n,\ga,s)$. 
Thinning out this system, and re-indexing, we can achieve 
$\ga+\la<x(n+1,\ga,s)$.

We can define $x(n+1,\ga,s\tau)<x(n,\ga,s)$ for $\tau<\la$ 
satisfying (1) and (2). 
Thinning and re-indexing, we can modify this system so that 
if we set $y_i=x(n+1,\ga,s|i)$ for $ i \leq n+1$, 
then 
$F(y_0,\dots,y_{n+1})=d(n+1)$ and 
$r_c(y_0,\dots,y_{n+1})=\xi(s,c)$ hold for every $s:n\to\la$, 
$c:\oo\to\mu$, i.e., the color and the rank do not depend on the last value. 

Repeating this, again thinning and re-indexing we finally get 
that the value of $r_c(y_0,\dots,y_{n+1})$  depends only on $c$, 
so it is a value $\xi(n+1,c)$, as claimed.

For the above function $d:\oo\to\mu$ we have that 
\[ 
      \xi(0,d)>\xi(1,d)>\cdots
\]
a contradiction. 
\qed 

\medskip
In order to handle scattered order types we represent them. 

\medskip
If $\a$ is an ordinal then let $H(\a)$ be the set of all 
$f:\a\to\{-1,0,1\}$ functions for which the set 
$D(f)=\{\b<\a:f(\b)\neq 0\}$ is finite. 
Order $H(\a)$ as follows. 
$f<f'$ iff $f(\b)<f'(\b)$ holds for the largest $\b$ with 
$f(\b)\neq f'(\b)$. 
This clearly orders $H(\a)$. 

\medskip
\noindent{\bf Lemma 2.} 
         {\em The order type of $\left(H(\a),<\right)$ is scattered.}

\medskip
\noindent{\bf Proof.}          
Assume that the mapping $q\to f_q$ is an order preserving 
injection for $q\in {\bf Q}$. 
Let $\b<\a$ be the {\em least} ordinal that occurs as the largest 
ordinal where $f_q$, $f_{q'}$ differ, for some $q<q'$. 
Now choose the rational numbers $q''$, $q'''$ with 
$q<q''<q'''<q'$. 
Then all four functions $f_{q}$,$f_{q'}$, $f_{q''}$,$f_{q'''}$ 
agree above $\b$, and some two at 
$\b$, too, a contradiction. 
\qed 
      
\medskip
\noindent{\bf Lemma 3.} 
         {\em Every scattered order type can be embedded into some 
         $\left(H(\a),<\right)$.}

\medskip
\noindent{\bf Proof.}          
Using Hausdorff's characterization it suffices to show that if 
some order types can be so represented then any well ordered and 
reverse well ordered sum of them can also be so represented. 
For this, it suffices to show that the antilexicographic products 
$H(\a)\times\b$ and $H(\a)\times\b^*$ can be embedded into 
$H(\a+\b)$. 
Indeed, if we map the pair $(f,\g)$ to the function $g$ 
which is $f$ restricted to $\a$ and in the interval $[\a,\a+\b)$ 
is everywhere zero except at $\a+\g$ where it is 1, then 
this is the required embedding for $H(\a)\times\b$. 
For the other case we use extensions that assume $-1$ at exactly one 
place. 
\qed 
      
\medskip
Given an $\a$-tree $\{x({\bf s}):{\bf s}\in \fs(\a)\}\subseteq\la^+$ 
we define an injection 
$\Phi:H(\a)\to H(\la^+)$ as follows. 
If $f\in H(\a)$, $D(f)=\{\b_0,\dots,\b_n\}$ in decreasing enumeration, 
then set 
$\g_j=x\left(\{\b_j,\dots,\b_0\}\right)$ for $0\leq j \leq n$. 
Now $\Phi(f)=g$ where $D(g)=\{\g_0,\dots,\g_n\}$ and $g(\g_j)=f(\g_j)$.

\medskip
\noindent{\bf Lemma 4.} 
         {\em This mapping $\Phi:H(\a)\to H(\la^+)$ is order preserving.}

\medskip
\noindent{\bf Proof.}      
Assume that $f$, $f'\in H(\a)$, 
$D(f)=\{\b_0,\dots,\b_n\}$, $D(f')=\{\b'_0,\dots,\b'_m\}$ 
in decreasing enumeration. 
Let $r$ be the largest index that 
for $i<r$ $\b_i=\b'_i$ and $f(\b_i)=f'(\b_i)$ hold.  
On some $\b$ we have $f(\b)<f'(\b)$ where either 
$\b=\b_r=\b'_r$ or 
$\b=\b_r\notin D(f')$ or 
$\b=\b'_r\notin D(f)$. 

Set $\g_j=x\left(\{\b_j,\dots,\b_0\}\right)$ for $j<r$ and 
$\g=x\left(\{\b,\b_{j-1},\dots,\b_0\}\right)$. 
Then the functions $\Phi(f)$ and $\Phi(f')$ agree above $\g$ and 
$\Phi(f)(\g)<\Phi(f')(\g)$ and we are done. 
\qed

\medskip
\noindent{\bf Theorem 2.} 
         {\em If $\phi$ is a scattered order type, $\mu$ is a cardinal, 
         then there exists a scattered order type $\psi$ such that 
         \[
         \psi \to [\phi]^{1}_{\mu,\aleph_0} 
         \]holds. }

\medskip
\noindent{\bf Proof.}          
By Lemmas  2., 3. it suffices to show that if $\a$ is an ordinal, 
$\mu$ a cardinal, then for some $\la$, the ordered set 
$\left(H(\la^+),<\right)$ has the property that for every coloring with 
$\mu$ colors there is a subset isomorphic to $\left(H(\a),<\right)$ which is 
colored with only countably many colors. 

Select $\la$ as in Theorem 1. 
Assume that $G:\left(H(\la^+),<\right)\to\mu$ is a coloring. 
Let $F$ be the following coloring of $\fs(\la^+)$. 
If ${\bf s}=s(0)s(1)\cdots s(n-1)$ is an element of it, 
let $F({\bf s})$ be the following function defined on 
$\{-1,1\}\times\cdots\{-1,1\}$. 
$F(i_0,\dots,i_{n-1})=G(f)$ where $f$ is the function with 
$D(f)={\bf s}$ and $f(s(j))=i_j$. 

Notice that this is a coloring with $\mu$ colors. 
By Theorem 1 there is an $\a$-tree 
$\{x({\bf s}):{\bf s}\in \fs(\a)\}$ such that 
\[ 
F\biggl( x\bigl(s(0)\bigr), x\bigl(s(0)s(1)\bigr), \cdots,
         x\bigl(s(0)s(1)\cdots s(n)\bigr) \biggr)=c(n).       
\]
holds for some function $c$.

If we now consider the corresponding mapping 
$\Phi:H(\a)\to H(\la^+)$ then it gives a subset of 
$\left(H(\la^+),<\right)$ isomorphic to 
$\left(H(\a),<\right)$ getting only $\mu$ colors. 
\qed

\vskip1cm
\hbox{
\vtop{\hbox{P\'eter Komj\'{a}th}
      \hbox{Department of Computer Science}
      \hbox{E\"otv\"os University}
      \hbox{Budapest, P.O.Box 120 } 
      \hbox{1518, Hungary}
      \hbox{e-mail:{\tt\ kope@cs.elte.hu}}}
      \qquad\qquad\quad
\vtop{\hbox{Saharon Shelah}
      \hbox{Institute of Mathematics,} 
      \hbox{Hebrew University,} 
      \hbox{Givat Ram, 91904,} 
      \hbox{Jerusalem, Israel}
      \hbox{e-mail:{\tt\ shelah@math.huji.ac.il}}}}


\medskip
\begin{thebibliography}{99}

\bibitem{erdhaj}P.~Erd\H os, A.~Hajnal: 
	On a classification of denumerable order types 
	and an application to the partition calculus, 	
   {\sl Fundamenta Mathematicae}, 	
   {\bf 51}(1962), 117--129. 

\bibitem{haus} F.~Hausdorff: 
	Grundz\"uge einer Theorie der Geordnete Mengen, 
	{\sl Math. Ann.}, {\bf 65}(1908), 435--505. 

\bibitem{rosen} Joseph G.~Rosenstein: 
	{\sl Linear orderings}, 
  	Academic Press, 1982. 



         
\end{thebibliography}
\end{document}